 \documentclass{gtart}
 \agtart
\input epsf
\def\relabelbox{%
  \hbox\bgroup%
}%
\def\endrelabelbox{%
\egroup%
}%
\def\relabel #1#2 {%
  \special{ps:/a {} def}%
  \smash{\rlap{#2}}%
}%
\def\adjustrelabel <#1,#2> #3#4 {%
  \special{ps:/a {} def}%
  \smash{\rlap{\kern #1 \raise #2\hbox{#4}}}%
}%
\def\extralabel <#1,#2> #3 {\smash{\rlap{\kern #1 \raise #2\hbox{#3}}}}%

\usepackage{geometry}                
\usepackage{graphicx}
\usepackage{amssymb}
\usepackage{epstopdf}
\usepackage{amssymb, amsmath, latexsym, amscd}
\usepackage{epsfig}
 \newtheorem*{theorem}{Theorem}
  \newtheorem*{corollary}{Corollary}
   \newtheorem*{lemma}{Lemma}

  \newcommand{\G}{\Gamma}
\newcommand{\bdry}{\partial}
 
\newcommand{\Zn}{Z_n}
\newcommand{\Zi}{Z_{\infty}}
\newcommand{\Mn}{M_n}
\newcommand{\Un}{U_n}
\newcommand{\Ui}{U_{\infty}}
\newcommand{\Gns}{\Gamma_{n,s}}
\newcommand{\Gn}{\Gamma_n}
\newcommand{\Gi}{\Gamma_\infty}

\newcommand{\me}{\mu\eta}

\title{Erratum to: Homology stability for outer automorphism groups of free groups}
\author{Allen Hatcher, Karen Vogtmann and Nathalie Wahl}
\addresses{Department of Mathematics, Cornell University, Ithaca, NY 14853\\ Department of Mathematics, Cornell University, Ithaca, NY 14853\\ Department of Mathematics, University of Chicago, Chicago, IL 60637}
\emails{hatcher@math.cornell.edu\\vogtmann@math.cornell.edu\\wahl@uchicago.edu }
 \begin{abstract} We correct the proof of Theorem 5 of the paper {\it Homology stability for outer automorphism groups of free groups}, by the first two authors.
 \end{abstract}
\begin{document}
\maketitle

In \cite{HatVog04} a proof was presented that the homology of certain groups $\Gamma_{n,s}$ is independent of both $n$ and $s$ for $n$ sufficiently large.  The groups $\Gamma_{n,s}$ include $Aut(F_n)$ ($s=1$) and $Out(F_n)$  ($s=0$).  In August of 2005 Nathalie Wahl discovered an error in the proof, and the purpose of this note is to fix that error. 

We assume the reader is familiar with \cite{HatVog04}, whose  notation and conventions we will use here without further comment.      
The error occurs in the first part of the  proof of Theorem 5, showing that the map $\beta_*\colon H_i(\Gamma_{n,s+2})\to H_i(\Gamma_{n+1,s})$ is injective for $n>>i$ and $s\ge 1$.   The argument used a diagram chase in the following diagram:
$$
\begin{CD}
H_i(\G_{n,s+2},\G_{n-1,s+4})@>>>H_{i-1}(\G_{n-1,s+4}) \\
@V{ \beta_*}VV@V{\beta_*}VV\\
H_i(\G_{n+1,s},\G_{n,s+2})@>>>H_{i-1}(\G_{n,s+2}) \\
\end{CD}
$$
It was asserted that the top horizontal and right vertical arrows were successive maps in the long exact sequence of the pair $(\G_{n,s+2},\G_{n-1,s+4})$, and hence their composition was the zero map, but in fact the group $\G_{n,s+2}$ in the lower right corner of the diagram is a different subgroup of  $\G_{n+1,s}$ from the $\G_{n,s+2}$ in the upper left corner, so that $\beta_*$ is not induced by the inclusion map of the pair.  
It is in fact true
that the composition is the zero map for $n$ sufficiently large, but a proof seems to require the results proved in this correction.

We correct the problem by giving a completely new proof of stability with respect to $s$ for $s \geq 1$, complementing the earlier proof of stability with respect to $n$. The new proof entirely avoids the diagram displayed above and instead focuses on the map $\mu\colon\Gns \to \G_{n,s+1}$.  In the range where $\alpha_*$ is an isomorphism, the relation $\alpha=\beta\mu^2$ derived in section 2 of \cite{HatVog04} shows that $\beta_*$ is an isomorphism if and only if $\mu_*$ is an isomorphism. Note that $\mu_*$ is always injective since it has a left inverse obtained by gluing a disk to one of the new boundary components, so we only need to prove surjectivity of $\mu_*$ in a stable range. The commutative diagram
$$
\begin{CD}
H_i(\G_{n,s})@>{\alpha_*}>>H_i(\G_{n+1,s})@>{\alpha_*}>>
\cdots@>{\alpha_*}>>H_i(\G_{n+k,s})\\
@VV{\mu_*}V @VV{\mu_*}V @. @VV{\mu_*}V\\
H_i(\G_{n,s+1})@>{\alpha_*}>> H_i(\G_{n+1,s+1})@>{\alpha_*}>>
\cdots@>{\alpha_*}>>H_i(\G_{n+k,s+1})
\end{CD}
$$
shows that if $\mu_*$ is an isomorphism for $n>>i$ then it will be an isomorphism in the same range that $\alpha$ is an isomorphism, namely $ n \geq 2i+2 $.  
It is thus not necessary to keep track of the precise stable range in the  arguments given in this correction.

  For $s\geq2$ let $\eta\colon \Gns\to \G_{n+1,s-1}$ be induced by gluing a copy of the 3-punctured sphere 
 $M_{0,3}$ to the first and last boundary components of $M_{n,s}$.  
The stabilization $\alpha$ is the composition $\eta\mu\colon\Gns\to\G_{n+1,s}$, and now we want to consider the opposite composition $\me\colon\Gns\to\G_{n+1,s}$. We will show that $\me$ is an isomorphism on $H_i$ for $ n>>i$, and hence $\mu_*\colon H_i(\G_{n,s-1}) \to H_i(\Gns)$ is surjective for $n>>i$ and $s\geq 2$. Since $ \alpha $ is a homology isomorphism for $n>>i$,  and, as we will see,  $\alpha$ commutes with $\me$, it will suffice to show that $\me$ is a homology isomorphism after passing to the direct limit under stabilization by $\alpha$. 
 This turns out to be much easier than showing $\me$ is a homology isomorphism before passing to the limit.

\begin{figure}[ht]
\centerline
{\relabelbox\small\epsfxsize 4truein
\epsfbox{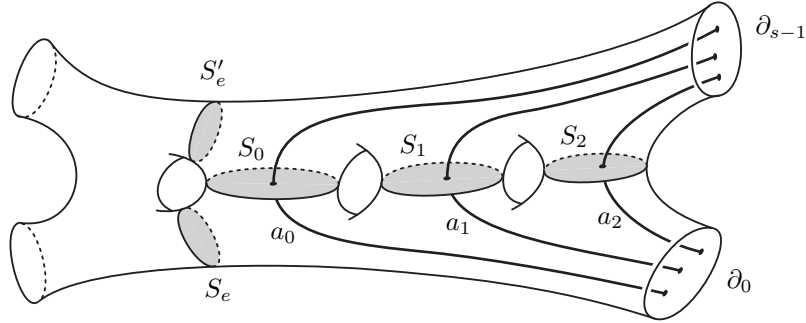}
\relabel{a}{$a_0$}
\relabel{b}{$a_1$}
\relabel{c}{$a_2$}
\relabel{bo}{$\partial_0$}
\relabel{bs}{$\partial_{s-1}$}
\relabel{S}{$S_0$}
\relabel{T}{$S_1$}
\relabel{U}{$S_2$}
\relabel{V}{$S_e$}
\relabel{W}{$S'_e$}
\endrelabelbox}
\caption{A 2-simplex in $Z_{n,s}$}
\label{simplex}
\end{figure}

To prove that $\me$ is a homology isomorphism in the limit we use a new simplicial complex $Z_{n,s}$, defined for $s\geq 2$.  To simplify the notation we will omit the subscript $s$ since it will be fixed throughout the proof, so we write $ Z_{n,s} $ as $ \Zn $ and $ M_{n,s} $ as $\Mn$. A vertex of $\Zn$ is an equivalence class of pairs $ (S,a) $ where $S$ is a non-separating sphere in $\Mn$ and $a$ is an embedded arc in $\Mn$ joining the first and last boundary spheres $\bdry_0$ and $\bdry_{s-1}$ and intersecting $S$ transversely in one point; we refer to $a$ as a {\it dual arc} for $S$.  The equivalence relation on such pairs $(S,a)$ is given by isotopy of $ S \cup a$ keeping the endpoints of $a$ in $\bdry\Mn$.  A set of $k+1$ vertices $(S_0,a_0),\cdots,(S_k,a_k)$ forms a $k$-simplex if $ S_i \cup a_i$ is disjoint from $ S_j \cup a_j$ for $ i \ne j $ and the spheres $S_i$ form a coconnected system (see Figure~\ref{simplex}). This implies that no two $S_i$'s or $a_i$'s are isotopic. Note that arcs can always be made disjoint by general position, so the disjointness condition really only involves intersections between different spheres and between spheres and arcs.

We note that there is an equivalent way of viewing a simplex of $\Zn$ in terms of enveloping sphere-pairs.  
 For a simplex $\{(S_0,a_0),\cdots,(S_k,a_k)\}$, take two parallel copies of each $S_i$, one on either side of $S_i$, and join each of these new spheres to a sphere parallel to either $\bdry_0$ or $\bdry_{s-1}$ by a tube following
the half of $a_i$ on the appropriate side of $S_i$.   
 This produces a pair of spheres $ S_e,S'_e$ separating $\Mn$ into two components, one of which contains the spheres $S_i$. The only boundary spheres contained in this component are $\bdry_0$ and $\bdry_{s-1}$, and splitting this component along the spheres $S_i$ produces two simply-connected pieces, one bounded by $S_e$, $\bdry_0$, and the $S_i$'s, the other bounded by $S'_e$, $\bdry_{s-1}$, and the $S_i$'s. Conversely, given a coconnected system $S_0,\cdots,S_k$ and two spheres $S_e,S'_e$ with the properties just listed, then there are dual arcs $a_i$ in the split-off submanifold such that $\{(S_0,a_0),\cdots,(S_k,a_k)\}$ is a simplex of $\Zn $, and these $a_i$'s are unique up to isotopy.

We wish to describe now an inclusion $\Zn \hookrightarrow Z_{n+1}$. This will be induced by an inclusion $\Mn \hookrightarrow M_{n+1}$. We have already used one such inclusion in the definition of the stabilization $\alpha\colon\G_{n,s}\to\G_{n+1,s}$ when we regarded $M_{n+1}$ as being obtained from $\Mn$ by attaching $M_{1,2}$ to $\bdry_0 $ along one boundary sphere of $M_{1,2}$. However, an alternative approach will make things a little clearer when dealing with the complexes $\Zn$. Here we build $M_{n+1}$ from $\Mn$ by  attaching $M_{1,1}$, identifying a disk in $\partial M_{1,1}$ with a disk in $ \bdry_0$. The first inclusion $\Mn \hookrightarrow M_{n+1}$ is then recovered by attaching a product $S^2 \times I$ to the new $\bdry_0$. Since attaching this product does not affect isotopy classes of diffeomorphisms modulo Dehn twists, the new inclusion $\Mn\hookrightarrow M_{n+1}$ gives the same $\alpha$ as the old one.  

\begin{figure}[ht] 
   \centering
   \includegraphics[width=4in]{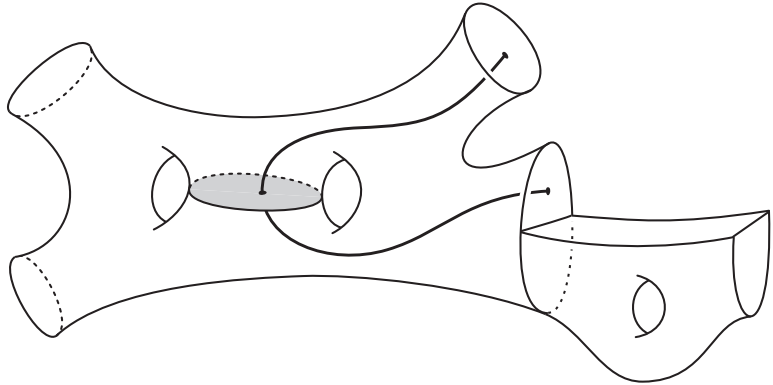} 
   \caption{Stabilization by $\alpha$ and the inclusion $\Zn\to Z_{n+1}$}
   \label{pairs}
\end{figure}

The new inclusion $\Mn \hookrightarrow M_{n+1}$ induces a map $\Zn \to Z_{n+1}$ since we may assume simplices of $\Zn$ are represented using pairs $(S_i,a_i)$ whose arcs $a_i$ are disjoint from the disk in $\bdry_0$ where $M_{1,1}$ is attached. This map $\Zn \to Z_{n+1}$ is injective by general properties of sphere systems (uniqueness of normal forms  \cite{Hat95}). Let $\Zi$ denote the direct limit of the complexes $\Zn$ under these inclusions $\Zn \hookrightarrow Z_{n+1}$.

\begin{lemma}  $Z_\infty$ is contractible.
\end{lemma}
\begin{proof} Given a map $g\colon S^k \to Z_\infty$, we wish to extend this to a map  $D^{k+1}\to Z_\infty$.  We may assume $g$ is simplicial with respect to some triangulation of $S^k$. This triangulation has finitely many simplices, so the image of $g$ lies in $\Zn$ for some $n$.  Let $M_n \subset M_{n+1} \subset \cdots \subset M_{n+k+1}$ be the alternate inclusions described above inducing  $\alpha$, and choose a non-separating sphere $T_i$ in each  $M_{n+i+1} - M_{n+i}$.

Triangulate $D^{k+1}$ by coning off the triangulation of $S^k$ to the centerpoint of $D^{k+1}$, a new vertex $v$. Define $g(v)$ to be $(T_0,b_v)$ where $b_v$ is any arc in $M_{n+1}$ dual to $T_0$.  Next, we extend $g$ over each interior edge $e$ of $D^{k+1}$ in the following way.  The endpoints of $e$ map to   $(T_0,b_v)$ and to another vertex  $(S_0,a_0)$.  Let $g$ send  the midpoint of $e$ to $(T_1,b_e)$, where $b_e$ is an arc in $M_{n+2}$ which is in the complement of the coconnected system $\{S_0,T_0\}$. Then $\{(S_0,a_0),(T_1,b_e)\}$ and $\{(T_1,b_e),(T_0,b_v)\}$ are edges of $Z_{n+2}$ so $g$ extends over $e$ by mapping its two halves to these two edges.

The extension of $g$ over simplices of $D^{k+1}$ of higher dimension proceeds in a similar fashion, by induction on the dimension of the simplices. Each $i$-simplex $\sigma$ of $D^{k+1}$ not contained in $S^k$ is the cone to $v$ of an $(i-1)$-simplex in $S^k$. The map $g$ sends this $(i-1)$-simplex to a possibly degenerate simplex $\{(S_0,a_0),\cdots,(S_{i-1},a_{i-1})\}$ in $\Zn$. The rest of the boundary of $\sigma$ is sent by induction to a subcomplex with additional vertices  $(T_j,b_\tau)$ for $0\leq j < i$, where $\tau$ ranges over the faces of $\sigma$ not in $S^k$. We send the barycenter of $\sigma$ to $(T_i,b_\sigma)$ where $b_\sigma$ is chosen in $M_{n+i+1}$ and in the complement of the coconnected system $\{S_0,\cdots,S_{i-1},T_0,\cdots,T_{i-1}\}$. We can then extend $g$ over $\sigma$ by coning off to its barycenter. This gives the induction step, and at the end of the induction we have extended $g$ over $ D^{k+1}$. Since $g$ and $k$ were arbitrary, this shows $\Zi$ is contractible.
\end{proof}

The natural action of $\Gn$  on $\Zn$ is transitive on simplices of each dimension, since splitting $M_n$ along the $k+1$ spheres of a $k$-simplex produces the manifold $M_{n-k-1}$ with $2k+2$ new punctures, each joined to $\bdry_0$ or $\bdry_{s-1}$ by an arc, and any two such configurations are diffeomorphic.  The action of $\Gn$ on $\Zn$ is compatible with the stabilization $\alpha$, in the sense that for each $ g \in \Gn $  the following diagram commutes:
$$\begin{CD}
\Zn @>g>>  Z_n \\
 @VVV  @VVV \\
Z_{n+1} @>{\alpha(g)}>>  Z_{n+1}
\end{CD}$$
Thus the direct limit group $\Gi$ acts on $\Zi$. This action is also transitive on $k$-simplices for each $k$. 

For the action of $\Gi$ on $\Zi$ the stabilizer of a simplex includes group elements that permute the vertices of the simplex, so to avoid this we consider $\Ui = \Delta(\Zi) $, the complex whose $k$-simplices are all the simplicial maps from the standard $k$-simplex to $\Zi$.  Note that $\Ui$ is the direct limit of the complexes $U_n=\Delta(Z_n)$. The homology of $\Ui$ is trivial since $\Zi$ has trivial homology. The action of $\Gi$ on $\Zi$ induces an action on $\Ui$. The quotient $\Ui/\Gi$ is contractible since it is the direct limit of the quotients $\Un/\Gn$ and these quotients are combinatorially the same as the quotients $W_n/\Gn$ in the proof of Theorem 4 of \cite{HatVog04}, which were $ (n-2)$-connected.

\begin{figure}[h] 
   \centering
   \includegraphics[width=4in]{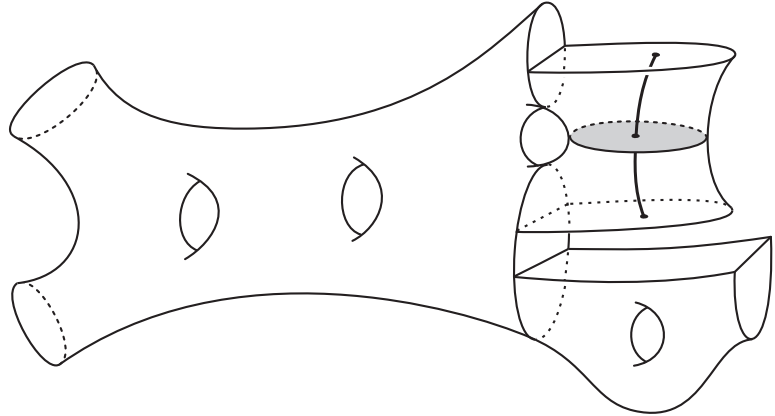} 
   \caption{Compatibility of $\me$ with $\alpha$}
   \label{stab2}
\end{figure}

The stabilizer of a vertex for the action of $\Gn$ on $\Un$ is a copy $\G_{n-1}'$ of $\G_{n-1}$ in $\Gn$, and the inclusion of this stabilizer is the map $\me$. The map $\me$ is induced by gluing a four-punctured 3-sphere to $\Mn$ by attaching two of its boundary spheres to $\bdry_0$ and $\bdry_{s-1}$.  As in the case of $\alpha$, there is an alternative description of $\me$ as being induced by gluing a  twice-punctured 3-sphere to $\Mn$ by attaching its boundary spheres to $\bdry_0$ and $\bdry_{s-1}$ along disks  (see Figure \ref{stab2}).

\noindent
This alternative description of $\me$ makes it clear that $\me$ commutes with $\alpha$, giving a commutative diagram:
$$
\begin{CD}
 \G_{n-1}' @>{\alpha}>> \Gn'@>{\alpha}>>\G_{n+1}' @>{\alpha}>> \G_{n+2}'@>{\alpha}>>\ldots\\
@VV{\me}V @VV{\me}V@VV{\me}V @VV{\me}V \ldots \\
\Gn @>{\alpha}>> \G_{n+1}@>{\alpha}>>\G_{n+2} @>{\alpha}>> \G_{n+3}@>{\alpha}>>\ldots\\
\end{CD}
$$
Thus in the limit action of $\Gi$ on $\Ui$ the inclusion of a vertex stabilizer is the direct limit map $ \me\colon\Gi' \to \Gi $.  Similarly, for stabilizers of higher dimensional simplices the inclusions of stabilizers are iterates of $\me$.

We can now use the equivariant homology spectral sequence arising from the action of $\Gi$  on $\Ui$ to prove:

\begin{theorem} The map $(\me)_*\colon H_i(\Gi')\to H_i(\Gi)$ is an isomorphism for each $ s\geq 2$.
\end{theorem}

As explained earlier, this implies:

\begin{corollary}
The map $\mu_*\colon H_i(\G_{n,s})\to H_i(\G_{n,s+1})$ is an isomorphism when $n\ge 2i+2$ and $s\ge 1$. 
\end{corollary}

\medskip\noindent{\it Proof of the theorem.}
 The proof proceeds by induction on $i$.  The equivariant homology spectral sequence has
 $$
 E^1_{p,q}=\bigoplus_{\sigma_p}H_q(stab(\sigma_p))\Rightarrow \widetilde H_{p+q}(\Ui)
 $$
 where $\{\sigma_p\}$ is a chosen set of orbit representatives for the $p$-simplices of $\Ui$.
 The differential $d^1\colon E^1_{0,i}\to E^1_{-1,i}$ is the map $\me\colon \Gi'\to \Gi$ we are interested in.  
 
The $j$-th row of the $E^1$ page of the spectral sequence is a chain complex computing the homology of the quotient $\Ui/\Gi$ with local coefficients in the system of groups $H_j(stab(\sigma_p))$.  
Each face of the boundary of $\sigma_p$ is equal to $h\sigma_{p-1}$ for some $\sigma_{p-1}$ and some element $h\in\Gi$, and the corresponding term of the $d^1$ map is the map $h_*\colon H_j(stab(\sigma_p))\to H_j(stab(\sigma_{p-1}))$ induced by conjugation by $h$.  For $j<i$ we may assume by induction that  the vertical maps in the commutative diagram below are isomorphisms.  
$$
\begin{CD}
H_j(stab(\sigma_p))@>h_*>>H_j(stab(\sigma_{p-1}))\\
@VVV @VVV\\
H_j(\Gi)@>h_*>>H_j(\Gi)\\
\end{CD}
$$
The lower $h_*$ in this diagram is the identity since it is induced by an inner automorphism of $\Gi$.  Thus   
the local coefficient system is trivial.  (Here we are following a line of reasoning that can be found in Section 7.4 of \cite{Iva87}.)

  Since the quotient has trivial homology, this shows that the entire $E^1$ page below the $i$-th row is zero.  The spectral sequence converges to 0 since $\Ui$ is contractible, and  the only differential with a chance of killing $E^1_{-1,i}$ is $d^1=\me$, proving that this map must be onto.  

To finish the induction we need to show that $\me$ is in fact an isomorphism on $H_i(\Gi')$. Since $\me$ is surjective on $H_i(\Gi')$, it is also surjective 
as a map $H_i(\G_{n-1,s}')\to H_i(\G_{n,s})$  
for large $n$ since $\alpha$ is an isomorphism on $H_i$ for large $n$.  Therefore $\mu_*\colon H_i(\G_{n,s-1})\to H_i(\Gns)$ is surjective for large $n$, and hence an isomorphism.  Since $\alpha=\eta\mu$ this implies that $\eta_*\colon H_i(\G_{n,s}')\to H_i(\G_{n+1,s-1})$ is an isomorphism for large $n$, so that $\me$ is also an isomorphism on $H_i(\G_{n,s}')$ for large $n$, hence on $H_i(\Gi')$ as well. \qed

\def\cprime{$'$}

\Addresses
\end{document}